\documentclass[11pt]{amsart}
\usepackage{amssymb}

\newcommand{\Lie}{\mathcal}
\newcommand{\nil}{\operatorname{nil}}
\newcommand{\prodsemi}{\rtimes}
\newcommand{\GL}{\operatorname{GL}}
\newcommand{\Ad}{\operatorname{Ad}}
\newcommand{\Aut}{\operatorname{Aut}}
\newcommand{\real}{\mathord{\mathbb R}}

\newcommand \ga{\closure{\Gamma^\alpha }}
\newcommand \G{\closure{\Ad G}}
\newcommand \Z[1]{\closure{\Ad_G #1}}
\newcommand \closure{\overline}
\newcommand \xc{\hat{X}}
\newcommand \zz[2]{\closure{\Ad_{#1} #2}}
\newcommand \z[1]{\closure{\Ad #1}}

\newcommand{\cover}{\widetilde}
\newcommand{\sigc}{{\cover{\sigma}}}
\newcommand{\gc}{{\cover{\ga}}}

\newcommand \stackarrow[1]{\setbox0\hbox{ \ $\scriptstyle#1$ \ }
 \setbox1\hbox{$\rightarrow$}\setbox2\hbox to
   \wd0{\rightarrowfill}\ht2=\ht1
 \mathrel{\mathop{\copy2}\limits^{\copy0}}}

\newcommand \pref[1]{\textnormal{(\ref{#1})}}
\renewcommand \see[1]{\textnormal{(}see~\textnormal{\ref{#1})}}

\renewcommand{\thmhead}[3]{%
 \thmname{#1}\thmnumber{ #2}\thmnote{ {\the\theoremnotefont#3}}}

\renewcommand{\swappedhead}[3]{%
 \textnormal{(\thmnumber{#2})}\thmname{ #1}\thmnote{
 {\the\theoremnotefont#3}}}
\swapnumbers

\newtheorem{thm}{Theorem}[section]

\newtheorem{lem}[thm]{Lemma}
\newtheorem{cor}[thm]{Corollary}
\newtheorem{prop}[thm]{Proposition}

\theoremstyle{definition}

\newtheorem{defn}[thm]{Definition}

\theoremstyle{remark}

\newtheorem{rem}[thm]{Remark}        
\newtheorem{ack}[thm]{Acknowledgment}

\newenvironment{thmref}{\thmrefer}{}
\newcommand{\thmrefer}[1]{\renewcommand\thethm
 {\protect\ref{#1}$'$}\addtocounter{thm}{-1}}

\newenvironment{thmcopy}{\thmcopier}{}
\newcommand{\thmcopier}[1]{\renewcommand\thethm
 {\protect\ref{#1}}\addtocounter{thm}{-1}}

\newbox\sectionS \setbox\sectionS\hbox{\S}

\newcommand{\refjour}{} %% make sure the names aren't already used
\newcommand{\refappear}{}
\newcommand{\refbook}{}

\long\def\refjour[#1]#2 #3: #4. #5.. #6 (#7) #8 \par
 {\bibitem{#2} #3, #4, \emph{#5} {\bf #6} (#7), #8.\par}
\long\def\refappear[#1]#2 #3: #4. #5 \par
 {\bibitem{#2} #3, #4, \emph{#5} (to appear).\par}
\long\def\refbook[#1]#2 #3: #4. #5, #6, #7 \par
 {\bibitem{#2} #3, ``#4,'' #5, #6, #7.\par}

\begin{document}

\title{Superrigid subgroups of solvable Lie groups}

\author{Dave Witte}

 \address{Department of Mathematics, Oklahoma State University, Stillwater, OK
74078}

\email{dwitte@math.okstate.edu}

\date{\today}

\thanks{\frenchspacing Submitted to
 \emph{Proceedings of the American Mathematical Society} (June 1996).}


\begin{abstract}
 Let $\Gamma$ be a discrete subgroup of a simply connected, solvable Lie group~$G$,
such that $\Ad_G\Gamma$ has the same Zariski closure as $\Ad G$. If $\alpha \colon
\Gamma \to \GL_n(\real)$ is any finite-dimensional representation of~$\Gamma $, we
show that $\alpha$ virtually extends to a continuous representation~$\sigma $
of~$G$. Furthermore, the image of~$\sigma$ is contained in the Zariski closure of
the image of~$\alpha $.
 When $\Gamma$ is not discrete, the same conclusions are true if we make the
additional assumption that the closure of $[\Gamma, \Gamma]$ is a finite-index
subgroup of $[G,G] \cap \Gamma$ (and $\Gamma$ is closed and $\alpha$ is
continuous).
 \end{abstract}

\maketitle

\section{Introduction}\label{intro}

The Margulis Superrigidity Theorem \cite[Thm.~VII.5.9, p.~230]{MargBook} concerns
lattices in semisimple Lie groups. The author \cite{Witte-super} (see
\cite[\S6]{Witte-super-S} for errata) proved an analogue of this fundamental theorem
for many lattices in Lie groups that are not semisimple, with most of his attention
devoted to solvable groups \see{super-latt}. We now prove that, in the case of
solvable groups, there is no need to restrict attention to lattices---exactly the
same result holds for any discrete subgroup.

\begin{defn} We say that a representation $\alpha \colon \Gamma\to
\GL_{n}(\real)$ \emph{virtually extends} to a representation of~$G$ if the restriction
of~$\alpha$ to some finite-index subgroup of~$\Gamma$ extends to a representation
of~$G$.
 \end{defn}

\begin{defn} \label{super-def}
 Let us say that a closed subgroup~$\Gamma$ of a simply connected, solvable Lie
group~$G$ is \emph{superrigid} in~$G$ if every finite-dimensional representation
$\alpha \colon \Gamma \to \GL_n(\real)$ virtually extends to a
representation~$\sigma$ of~$G$, such that the image $G^\sigma$ is contained in the
Zariski closure of~$\Gamma ^\alpha $.
 \end{defn}

\begin{rem} \label{fin-ind-super}
 By considering induced representations, we see that every finite-index
subgroup of a superrigid subgroup is superrigid.
 \end{rem}

\begin{thmref}{Zar-def}
 \begin{defn}
 If $X$ is a subgroup of a Lie group~$G$, let\/ $\Z{X}$ denote the (almost) Zariski
closure of\/ $\Ad_G X$ in the real algebraic group\/ $\Aut \Lie G$, where $\Lie G$ is
the Lie algebra of~$G$. 
 \end{defn}
 \end{thmref}

\begin{thm}[{\see{super-discrete}}] \label{weak-discrete}
 Let\/ $\Gamma$ be a discrete subgroup of a simply connected, solvable Lie group~$G$,
such that\/ $\Z{\Gamma} = \G$. Then\/ $\Gamma$ is superrigid in~$G$.
 \end{thm}

For subgroups that are not discrete, an additional hypothesis is required.

\begin{thm}[{\see{super-closed}}] \label{easy-super}
 Let\/ $\Gamma$ be a closed subgroup of a simply connected, solvable Lie group~$G$,
such that\/ $\Z{\Gamma} = \G$. If the closure of\/ $[\Gamma, \Gamma]$ is a
finite-index subgroup of\/ $[G,G] \cap \Gamma $, then\/ $\Gamma$ is superrigid in~$G$.
 \end{thm}

If $\Gamma$ is superrigid in~$G$, then $\Gamma$ is also superrigid in any semidirect
product $A \prodsemi G$. The following theorem shows that combining this observation
with Thm.~\ref{easy-super} is (virtually) the only way to construct a superrigid
subgroup.

\begin{thmcopy}{super-converse}
 \begin{thm}
 A closed subgroup\/~$\Gamma$ of a simply connected, solvable Lie group~$G$ is
superrigid
 iff\/
 $\Gamma$ has a finite-index subgroup\/~$\Gamma'$, such that
 \begin{enumerate}
 \item
 there is a semidirect-product decomposition $G = A \prodsemi B$, with\/ $\Gamma'
\subset B$ and\/
 $\zz{B}{\Gamma'} = \z{B}$, and
 \item the closure of\/ $[\Gamma', \Gamma']$ is a
finite-index subgroup of\/ $[B,B] \cap \Gamma' $.
 \end{enumerate}
 \end{thm}
 \end{thmcopy}

If $\Gamma$ is discrete, then Hypothesis~\pref{must-fin-ind} of the theorem follows
from~\pref{must-semi-prod} (see~\ref{derived-cocpct}\pref{fin-ind-derived}), so we
have the following corollary.

\begin{cor} \label{cnvs-discrete}
 A discrete subgroup\/~$\Gamma$ of a simply connected, solvable Lie group~$G$ is
superrigid
 iff
 there is a semidirect-product decomposition $G = A \prodsemi B$, such that $B$
contains a finite-index subgroup\/~$\Gamma'$ of\/~$\Gamma $, and\/
 $\zz{B}{\Gamma'} = \z{B}$. \qed
 \end{cor}

\begin{ack}
 I am grateful to A.~Magid for asking the astute question that directed me to this
extension of my earlier work. This research was partially supported by a grant from
the National Science Foundation.
 \end{ack}

\section{Preliminaries}\label{prelims}

\begin{defn}[{\cite[Defn.~3.2]{Witte-super}}]
 A subgroup~$A$ of $\GL_n(\real )$ is said to be \emph{almost Zariski closed} if
there is a Zariski closed subgroup~$B$ of $\GL_n(\real )$, such that $B^\circ
\subset A \subset B$, where $B^\circ$ is the identity component of~$B$ in the
topology of $\GL_n(\real)$ as a $C^\infty$~manifold (not the Zariski
topology). There is little difference between being Zariski closed and almost
Zariski closed, because $B^\circ$ always has finite index in~$B$.
 \end{defn}

\begin{defn}[{\cite[Defn.~3.6]{Witte-super}}] \label{Zar-def}
 The \emph{almost-Zariski closure}~$\closure{A}$ of a subgroup~$A$ of
$\GL_n(\real)$ is the unique smallest almost-Zariski closed subgroup
that contains~$A$.
 In particular, if $A$ is a subgroup of a Lie group~$G$, we use
$\Z{A}$ to denote the almost-Zariski closure of $\Ad_G A$ in the real
algebraic group $\Aut(\Lie G)$, where $\Lie G$ is the Lie algebra
of~$G$.
 \end{defn}

\begin{defn}[({\cite[\S5]{Witte-super},
cf.~\cite[p.~6]{Fried-Goldman}})]
 \label{synd-def}
 Let $\Gamma$ be a closed subgroup of a Lie group~$G$. A {\it syndetic hull}
of~$\Gamma$ is a connected subgroup~$B$ that contains~$\Gamma$, such that
$B/\Gamma$ is compact.
 \end{defn}

\begin{prop}[{\cite[Prop.~5.9, p.~168]{Witte-super}}] \label{synd-dense}
 Let\/ $\Gamma$ be a closed subgroup of a connected, solvable Lie group~$G$, such
that\/ $\Z{\Gamma}$ contains a maximal compact torus of\/ $\G$. Then\/ $\Gamma$ has a
syndetic hull in~$G$. \qed
 \end{prop}

\begin{lem}[{\cite[Prop.~5.10 and Cor.~5.11, p.~168]{Witte-super}}]
\label{derived-cocpct}
 Suppose\/ $\Gamma$ is a closed subgroup of a simply connected, solvable Lie
group~$G$, such that\/ $\Z{\Gamma} = \G$. Then the closure of\/ $[\Gamma, \Gamma]$
contains a cocompact subgroup of\/ $[G,G] \cap \Gamma$. Furthermore:
 \begin{enumerate}
 \item \label{fin-ind-derived}
 If\/ $\Gamma$ is discrete, then\/ $[\Gamma, \Gamma]$ is a finite-index
subgroup of\/ $[G,G] \cap \Gamma $.
 \item \label{cont-derived}
 If\/ $\Gamma$ is connected, then\/ $[\Gamma, \Gamma] = [G,G]$. \qed
 \end{enumerate}
 \end{lem}

\begin{defn}
 A real-valued function $f$ on a subset~$G$ of $\GL_m(\real )$ is \emph{regular} if
there are polynomials~$p$ and~$q$ in $n^2 + 1$ variables, such that, for all $g \in
G$, we have
 $$
 q(g_{ij}, \det g ^{-1} ) \not= 0
\qquad \textnormal{and} \qquad
 f(g) = \frac{ p(g_{ij}, \det g ^{-1} )}{q(g_{ij}, \det g ^{-1} )} .
 $$ 
 A function $f \colon G \to \GL_n(\real)$ is \emph{regular} if each matrix entry of
$f(g)$ is a regular function of~$g$.
 \end{defn}

\begin{lem}[{\cite[Prop.~3.16, p.~158]{Witte-super}}] \label{image-closed}
 Suppose $G$ is an almost-Zariski closed subgroup of\/ $\GL_m(\real)$, and $\phi
\colon G \to \GL_n(\real)$ is a regular homomorphism. Then $G^\phi$ is almost-Zariski
closed. \qed
 \end{lem}

\begin{lem}[{\cite[Lem.~3.17, p.~159]{Witte-super}}] \label{prod-closed}
 Suppose $A$~and~$B$ are almost-Zariski closed subgroups of\/ $\GL_n(\real)$, such
that $AB$ is a subgroup. Then $AB$ is almost Zariski closed. \qed
 \end{lem}

\begin{lem}[{\cite[Cor.~2.19, p.~154]{Witte-super}}] \label{quot-sc}
 If $N$ is a connected, normal subgroup of a simply connected, solvable Lie
group~$G$, then $G/N$ is simply connected. Therefore, $G/N$ has no nontrivial compact
subgroups. \qed
 \end{lem}

\begin{lem}[(Ado-Iwasawa)] 
\label{Ado-Thm}
 If $G$ is a simply connected, solvable Lie group, then $G$ has a faithful,
finite-dimensional representation~$\alpha$, such that $G^\alpha$ is closed.
 \end{lem}

\begin{proof}
 The usual statement of the Ado-Iwasawa Theorem \cite[Thm.~XVIII.3.1,
p.~219]{Hochschild} states that $G$ has a faithful, finite-dimensional
representation~$\sigma $, such that $(\nil G)^\sigma$ is unipotent. Hence $(\nil
G)^\sigma$ is closed \see{unip-closed}. For the same reason, because $G/\nil G$ is
abelian, there is faithful representation~$\tau$ of $G/\nil G$, whose image is 
closed. The direct sum of~$\sigma$ and~$\tau$ is the desired representation~$\alpha $.
 \end{proof}

\begin{lem}[{\cite[Lem.~5.6, p.~166]{Witte-super}}] \label{inv-conn}
 Let $G$ be a connected, solvable Lie group, and let~$A$ be an almost-Zariski closed
subgroup of\/ $\GL_n(\real)$. If $\rho \colon G \to \GL_n(\real)$ is a continuous
homomorphism, such that $A$ contains a maximal compact torus of~$\closure{G^\rho}$,
then the inverse image $\rho ^{-1} (A)$ is a connected subgroup of~$G$. \qed
 \end{lem}

\begin{lem} \label{ker-conn}
 Let $G$ and $H$ be simply connected, solvable Lie groups. If $\sigma \colon G \to
H$ is a continuous homomorphism, then the kernel of~$\sigma$ is connected.
 \end{lem}

\begin{proof}
 Let $K = \ker\sigma $. Because every connected subgroup of a simply connected,
solvable Lie group is simply connected \cite[Thm.~XII.2.2, p.~137]{Hochschild}, we
know $G^\sigma$ is simply connected. In other words, $G/K$ is simply
connected. Because $G$ is connected and solvable, this implies that $K$ is connected
\cite[Lem.~2.17, p.~154]{Witte-super}.
 \end{proof}

\begin{lem}[{\cite[Lem.~3.20, p.~159]{Witte-super}}] \label{unip-closed}
 Every connected, unipotent Lie subgroup of\/ $\GL_n(\real)$ is Zariski closed. \qed
 \end{lem}

\begin{lem}[{(cf.~\cite[Cor.~I.5.3.7, p.~47]{Bourbaki})}] \label{nilrad-unip}
 Let $A$ be a connected, nilpotent, normal subgroup of a Lie group~$G$. Then\/
$\Z{\nil G}$ is unipotent. \qed
 \end{lem}

The following well-known result is a consequence of the fact that maximal compact
tori are conjugate (see \cite[Prop.~3.10(2), p.~126]{Platonov-Rapinchuk} and
\cite[Cor.~102.9.1, p.~293]{Zelobenko}).

\begin{lem} \label{normal-tor}
 Let $K$ and $G$ be almost-Zariski closed, connected subgroups of
$\GL_n(\real)$.
 If $K$ is a normal subgroup of~$G$, then every maximal compact torus of~$G$ contains
a maximal compact torus of~$K$. \qed
 \end{lem}

\section{The theorems and their proofs}\label{main-thm-sect}

The following are two versions of the main theorem of \cite{Witte-super}. In this
section, we eliminate the cocompactness assumption in these results \see{super-rig}
and prove a converse \see{super-converse}.

\begin{thm}[{\cite[Cor.~6.6, p.~175]{Witte-super}}] \label{cocpct-extend}
 Let\/ $\Gamma$ be a closed, cocompact subgroup of a simply connected, solvable Lie
group~$G$, and assume\/ $\Z{\Gamma} = \G$. Suppose $\alpha \colon \Gamma \to
\GL_n(\real)$ is a representation of\/~$\Gamma$. If\/ $\ga$ and\/ $Z(\ga)$ are
connected, and\/ $([G,G] \cap \Gamma)^\alpha$ is unipotent, then $\alpha$ extends to
a homomorphism from~$G$ to\/~$\ga$. \qed
 \end{thm}

\begin{cor}[{\cite[Cor.~6.8, p.~176]{Witte-super}}] \label{super-latt}
 Let\/ $\Gamma$ be a lattice subgroup of a simply connected, solvable Lie
group~$G$, such that\/ $\Z{\Gamma} = \G$. Then\/ $\Gamma$ is superrigid in~$G$. \qed
 \end{cor}

\begin{rem}
 If one is interested only in \pref{weak-discrete}, not the more
precise results such as Thm~\ref{super-rig} and its corollaries, a
direct proof can be obtained from the proof of~\pref{super-rig} by appealing to
Cor.~\ref{super-latt} instead of Thm.~\ref{cocpct-extend}. A short, fairly easy proof
of Cor.~\ref{super-latt} appears in \cite[\S2]{Witte-super-S}.
 \end{rem}

\begin{thm} \label{super-rig}
 Let $\Gamma$ be a closed subgroup of a simply connected, solvable Lie group~$G$, and
assume $\Z{\Gamma} = \G$. Suppose $\alpha \colon \Gamma \to
\GL_n(\real)$ is a representation of\/~$\Gamma$. If\/ $\ga$ and\/ $Z(\ga)$ are
connected, and\/ $([G,G] \cap \Gamma)^\alpha$ is unipotent, then $\alpha$ extends to
a homomorphism from~$G$ to\/~$\ga$.
 \end{thm}

\begin{proof} Because $\Z{\Gamma} = \G$, we know that $\Gamma$ has a syndetic
hull~$B$ in~$G$ \see{synd-dense}. Then $\alpha$ extends to a homomorphism $\sigma 
\colon B \to \ga$ \see{cocpct-extend}. From
Lem.~\ref{derived-cocpct}\pref{cont-derived}, we see that 
 $([G,G] \cap B)^\sigma = [B,B]^\sigma \subset [\ga, \ga]$ 
 is unipotent, so there is no harm in replacing~$\Gamma$ with~$B$. Thus, we may assume
$\Gamma$ is connected.

We may assume $G \subset \GL_n(\real)$ \see{Ado-Thm}, so there is no harm in speaking
of the almost-Zariski closure of~$G$ or of its subgroups. Let $\pi \colon
\closure{G} \times \ga \to \closure{G}$ be the projection onto the first factor.

Let $X \subset \Gamma \times \ga$ be the graph of~$\alpha $. It suffices to find
a closed subgroup~$Y$ of $G \times \ga$ such that
 $X \subset Y$,
 $Y \cap \ga = e$, and
 $Y^\pi =G$.
 For then $Y$ is the graph of a homomorphism $G \to \ga$ that extends~$\alpha $, as
desired.

Because $\closure{X}^\pi \supset X^\pi =\Gamma $, and $\Ad_G \closure{X}^\pi$ is
almost Zariski closed \see{image-closed}, we have
 $\Ad_G \closure{X}^\pi \supset \Z{\Gamma} = \G$,
 so $G \subset \closure{X}^\pi Z(\closure{G})$.
 Thus, letting $\xc = \closure{X} \, Z(\closure{G})^\circ$,
 we have $G \subset \xc^\pi $.

Let $Y$ be a connected subgroup of~$\xc$ that contains~$X$ and is maximal, subject to
the conditions that $Y \cap \ga = e$ and $Y^\pi \subset G$. We claim that $Y^\pi =G$
(which will complete the proof). If not, then, because $G \subset \xc^\pi $, there is
a one-parameter subgroup~$A$ of~$\xc$ such that $A^\pi$ is contained in~$G$, but is
not contained in~$Y^\pi $. From Lem.~\ref{derived-cocpct}\pref{cont-derived}, we know
that
 $[ \closure{X}, \closure{X}] = [X,X] \subset X \subset Y$,
 so $Y$ is normal in~$\xc$; hence $AY$ is a subgroup. Furthermore, because $G$ is
simply connected and $Y^\pi$ is connected, we know that $G/Y^\pi$ has no nontrivial
compact subgroups \see{quot-sc}. Because $A^\pi \not \subset Y^\pi $, this implies
that $a^\pi \not\in Y^\pi $, for all $a \in A- \{e\}$; therefore $(AY) \cap \ga = Y
\cap \ga$ is trivial. This contradicts the maximality of~$Y$.
 \end{proof}

\begin{cor}[{(cf.~pf.~of \cite[Cor.~6.7, p.~176]{Witte-super})}]
\label{super-discrete}
 Let\/ $\Gamma$ be a discrete subgroup of a simply connected, solvable Lie group~$G$,
and assume\/ $\Z{\Gamma} = \G$. Suppose $\alpha \colon \Gamma \to \GL_n(\real)$ is a
homomorphism, such that\/ $\ga$ is connected. Then there is a homomorphism $\sigma
\colon G \to \ga$ and a finite subgroup~$Z_0$ of\/ $Z(\ga)$ such that $\gamma
^{\sigma} \in \gamma ^\alpha Z_0$, for every $\gamma \in \Gamma$. In particular,
$\sigma$ virtually extends~$\alpha $. \qed
 \end{cor}

To obtain the same conclusion for subgroups that are not discrete, we add the
conclusion of Lem.~\ref{derived-cocpct}\pref{fin-ind-derived} as a hypothesis.

\begin{cor}[{(cf.~pf.~of \cite[Cor.~6.7, p.~176]{Witte-super})}] \label{super-closed}
 Let\/ $\Gamma$ be a closed subgroup of a simply connected, solvable Lie group~$G$,
such that the closure of\/ $[\Gamma, \Gamma]$ is a finite-index subgroup of\/ $[G,G]
\cap \Gamma$, and assume\/ $\Z{\Gamma} = \G$. Suppose $\alpha \colon \Gamma \to
\GL_n(\real)$ is a homomorphism, such that\/ $\ga$ is connected. Then there is a
homomorphism $\sigma \colon G \to \ga$ and a finite subgroup~$Z_0$ of $Z(\ga)$ such
that $\gamma ^{\sigma} \in \gamma ^\alpha Z_0$, for every $\gamma \in \Gamma$. In
particular, $\sigma$ virtually extends~$\alpha $. \qed
 \end{cor}

Without the additional hypothesis, we obtain only the following weaker result, by
using the fact that, although the closure of $[\Gamma, \Gamma]$ may not have finite
index in $[G,G] \cap \Gamma $, it is cocompact \see{derived-cocpct}.

\begin{cor}[{(cf.~pf.~of \cite[Cor.~6.7, p.~176]{Witte-super})}]
 Let $\Gamma$ be a closed subgroup of a simply connected, solvable Lie group~$G$, and
assume $\Z{\Gamma} = \G$. Suppose $\alpha \colon \Gamma \to \GL_n(\real)$ is a
homomorphism, such that $\ga$ is connected. Then there is a homomorphism $\sigma
\colon G \to \ga$ and a compact subgroup~$Z_0$ of $Z(\ga)$ such that $\gamma ^{\sigma}
\in \gamma ^\alpha Z_0$, for every $\gamma \in \Gamma$. \qed
 \end{cor}

\begin{thm} \label{super-converse}
 A closed subgroup\/~$\Gamma$ of a simply connected, solvable Lie group~$G$ is
superrigid iff\/
 $\Gamma$ has a finite-index subgroup~$\Gamma'$, such that
 \begin{enumerate}
 \item \label{must-semi-prod}
 there is a semidirect-product decomposition $G = A \prodsemi B$, with\/ $\Gamma'
\subset B$ and\/
 $\zz{B}{\Gamma'} = \z{B}$, and
 \item \label{must-fin-ind}
 the closure of\/ $[\Gamma', \Gamma']$ is a finite-index subgroup of\/ $[B,B] \cap
\Gamma' $.
 \end{enumerate}
 \end{thm}

\begin{proof} ($\Rightarrow$)
\pref{must-semi-prod}
 There is a faithful representation $\alpha \colon G \to \GL_n(\real)$, such that
$G^\alpha$ is closed \see{Ado-Thm}. (So $\Gamma ^\alpha$ is also closed.) Because
$\Gamma$ is superrigid, we know that, after replacing~$\Gamma$ by a finite-index
subgroup, the homomorphism $\alpha |_\Gamma$ extends to a homomorphism $\sigma 
\colon G \to \ga$. Because $G$ is simply connected, we may lift~$\sigma$ to a
homomorphism $\sigc \colon G \to \gc$, where $\gc$ is the universal cover of~$\ga$.
By combining Lem.~\ref{image-closed} with the fact that $\ga = \closure{G^\sigma}$,
we see that
 $$
 \zz{G^{\sigma}}{\Gamma ^\alpha} 
 = \Ad_{G^\sigma} \ga
 = \Ad_{G^\sigma} \closure{G^\sigma}
 = \z{G^\sigma}.$$
 Therefore, by passing to the covering group~$G^\sigc$, we have
 $\zz{G^\sigc}{\Gamma ^\sigc} = \z{G^\sigc}$.
 Thus, letting $K = \ker \sigc$, we see that $\Z{K\Gamma} = \G$. Furthermore, $K$ is
connected \see{ker-conn}, we have $K \cap \Gamma =e$, because $\sigma|_\Gamma =
\alpha$ is faithful, and we know that $K\Gamma$ is closed, because $\Gamma ^\alpha$
is closed. So Lem.~\ref{semi-prod} below implies that there is a semidirect-product
decomposition $G = A \prodsemi B$ with $K \subset A$ and $\Gamma \subset B$. Because
$\Z{A\Gamma } \supset \Z{K\Gamma} = \G$, we see that $\zz{B}{\Gamma} = \z{B}$.

 \pref{must-fin-ind}
 Because $\Gamma'$ is a finite-index subgroup of~$\Gamma $, we know that $\Gamma'$
is a superrigid subgroup of~$G$ \see{fin-ind-super}. Thus, by replacing~$\Gamma $
with~$\Gamma'$, we may assume $\Gamma \subset B$.
 Let $D$ be the closure of $[\Gamma, \Gamma]$. Because $\Gamma /D$ is abelian, there
is a faithful homomorphism $\alpha \colon \Gamma /D \to T$, for some compact torus
$T \subset \GL_n(\real)$. We may think of~$\alpha$ as a representation of~$\Gamma$
that is trivial on~$D$. Then, because $\Gamma$ is superrigid, we know that $\alpha$
virtually extends to a homomorphism $\sigma \colon G \to T$. Of course $[G,G]^\sigma
=e$, because $T$ is abelian. Thus, because $\alpha$ agrees with~$\sigma$ on a
finite-index subgroup of~$\Gamma $, we conclude that~$D$, the kernel of~$\alpha $,
contains a finite-index subgroup of $[G,G] \cap \Gamma $.
 \end{proof}

\begin{lem} \label{semi-prod}
 Let $K$ and\/ $\Gamma$ be closed subgroups of a simply connected, solvable Lie
group~$G$, such that $K$ is connected and normal in~$G$, $K \cap \Gamma = e$,
$K\Gamma$ is closed, and\/ $\Z{K\Gamma} = \G$. Then there is a semidirect-product
decomposition $G = A \prodsemi B$ with $K \subset A$ and $\Gamma \subset B$.
 \end{lem}

\begin{proof}
 By induction on the derived length of~$K$, we may assume $K$ is abelian, so $\Ad_G
K$ is unipotent \see{nilrad-unip}. Then, because $\Z{K\Gamma} = \G$, we know that
$\Z{\Gamma}$ contains a maximal compact torus of $\G$. So Lem.~\ref{synd-dense}
implies that $\Gamma$ has a syndetic hull~$X$ in~$G$. Thus, by replacing $\Gamma$
with~$X$, we may assume $\Gamma$ is connected.

Let $H = \{ g \in G \mid \Ad_G g \in \Z{\Gamma} \}$. Because $\Z{\Gamma}$ contains a
maximal compact torus of $\G$, we know that $H$ is connected \see{inv-conn}. For much
the same reason, $K \cap H$ is connected \see{normal-tor}.

We claim that $G = K H$. Given $g \in G$, because $\G = \Z{K\Gamma}$, we see, from
Lem.~\ref{prod-closed}, that $\Ad_G g = k' h'$, for some $k' \in \Z{K}$ and $h' \in
\Z{\Gamma}$. But, because $\Ad_G K$ is unipotent, we have $\Ad_G K = \Z{K}$
\see{unip-closed}, so there is some $k \in K$ with $\Ad_G k = k'$. Then $\Ad_G(k^{-1}
g) = h' \in \Z{\Gamma}$, so, from the definition of~$H$, we have $k^{-1} g \in H$.

From the preceding paragraph, we see that, by replacing $K$ with $K \cap H$, we may
assume $G = H$. Therefore, $\Z{\Gamma} = \G$, so $[G,G] \subset \Gamma$
(see~\ref{derived-cocpct}\pref{cont-derived}). Thus, there is no harm in modding out
$[G,G]$, so we may assume $G$ is abelian. In this case, it is easy to obtain the
desired conclusion.
 \end{proof}

\end{document}